%% file: GCASfinal3.tex
\documentclass{article}

\usepackage{amssymb,amsmath}
\usepackage{latexsym}
\usepackage[all]{xy}

\input{amsmacros}

\def\mg{M_g}
\def\mgbar{\overline{M}_g}

\def\stg{\mathcal{M}_g}
\def\stgbar{\overline{\mathcal{M}}_g}

\begin{document}

\title{General curves on algebraic surfaces}
\author{ E. SERNESI \footnote{The author is member of GNSAGA-INDAM.
Work supported by the 2008 PRIN project  \emph{Geometria delle variet\`a algebriche e loro spazi di moduli}.}}
\date{}
\maketitle

\abstract{We give upper bounds on the genus of a curve with
general moduli assuming that it can be embedded in a projective
nonsingular surface $Y$ so that $\dim(|C|) > 0$. We find such
bounds for all types of surfaces of intermediate Kodaira dimension
and, under mild restrictions,  for surfaces of general type whose
minimal model $Z$ satisfies the Castelnuovo inequality $K_Z^2 \ge
3\chi(\O_Z) - 10$. In this last case we obtain $g \le 19$. In the
other  cases considered the bounds are lower.}

\section*{Introduction}

\textbf{Notation:} $\mgbar$ is the coarse moduli  space of
 stable   curves of genus $g$,  and $\mg
\subset \mgbar$ is the open subset of smooth curves.  The
corresponding \emph{stacks} are denoted by $\stgbar$ and $\stg$
respectively.
 We work over $\mathbb{C}$.

\separation

  In this paper we  study the conditions   imposed on $g$ by assuming that
  a general projective nonsingular curve $C$ of genus $g$
can be embedded in some non-ruled algebraic surface $Y$ so that
$\dim(|C|) > 0$ on $Y$. We expect upper bounds on $g$ depending on
the numerical characters of $Y$, because the assumption made
implies that $\mgbar$  is uniruled.  In fact
 blowing-up the base locus of a linear pencil extracted from $|C|$
one constructs a non-isotrivial fibration in curves of genus $g$
\begin{equation}\label{E:fib}
 f:X\longrightarrow
\P^1
\end{equation}
 containing $C$ among its fibres, which defines a non-constant morphism
 $\psi_f:\P^1\longrightarrow \mgbar$ whose image contains $[C]$.
 Conversely every such
fibration determines a pair $C\subset Y$ as above (in several
ways).  In general $f$ is not semistable, and we cannot apply
semistable reduction because rationality of the base would be lost
after the process.  We remedy to this drawback by   working
directly with an arbitrary non-isotrivial fibration over $\P^1$
and by studying its deformation theory. The deformation functor of
a fibration (\ref{E:fib}) is controlled by the sheaf
$Ext^1_f(\Omega^1_{X/\P^1},\O_X)$ on $\P^1$. A necessary condition
for (\ref{E:fib})  to deform enough to contain a general curve of
genus $g$ as a fibre is that this sheaf is globally generated
(Proposition \ref{C:free1} and Theorem \ref{T:unirfree}), and in
this case we call $f$ a \emph{free fibration}. We prove
(Proposition \ref{P:def2}) that if (\ref{E:fib}) is free then the
following inequality is satisfied:
\begin{equation}\label{E:ineq1}
    11\chi(O_X)- 2 K_X^2 \ge 5(g-1)+2
\end{equation}
 If we know that the  fibration
(\ref{E:fib}) is free and is obtained by blowing-up the base locus
of a linear pencil extracted from a linear system $|C|$ of
positive dimension   on a surface $Y$, then the above  inequality
can be refined
  as follows:
  \begin{equation}\label{E:ineq2}
    10\chi(O_Y)-2K_Y^2 \ge 4(g-1)-C^2-h^0(K_Y-C)
\end{equation}
(Theorem \ref{T:def1}). This inequality provides a necessary
condition for the possibility of embedding a general curve of
genus $g$ in a surface $Y$ with the given numerical characters so
that dim$(|C|)>0$.

In the second part of the paper we apply inequality
  (\ref{E:ineq2})  to give upper bounds on $g$
depending on the numerical characters of the surface $Y$, for such
an embedding being possible.  The case when $Y$ is rational has
been considered classically by   B. Segre \cite{bS29}. His work
has been  reconsidered and surveyed by Verra in his recent
\cite{aV11}. It gives the bound $g \le 10$ for a general curve of
genus $g$   moving in a non-trivial linear system of plane curves,
provided the system is regular. It has been conjectured that this
condition is always satisfied. We do not consider this case in the
present work, referring the reader to \cite{aV11} for a detailed
discussion and to \cite{eS12} for a general survey.

 Recent work by Verra and Bruno-Verra \cite{aV05,BV05},
plus experimental evidence in the range $12 \le g \le 15$, seem to
indicate that if $\mgbar$ is uniruled for some $g \ge 17$ then a
general curve of genus $g$ should move in a non-trivial linear
system on a regular surface \emph{of general type.}   We give
evidence for this expectation  by proving a series of results.
For elliptic surfaces we have:

\begin{theorem}\label{T:intr2}
Let $Y$ be a non-ruled and non-rational   elliptic surface, $\pi:Y
\to B$ the elliptic fibration onto a nonsingular connected curve
$B$.  Assume that there exists    a general nonsingular connected
curve  $C \subset Y$  of  genus $g\ge 3$ such that dim$(|C|) \ge
1$.  Then $B=\P^1$ and  $g \le 16$.
\end{theorem}

  The following result covers the  case of curves on other surfaces of intermediate
Kodaira dimension.

\begin{theorem}\label{T:intr3}
Let $Y$ be a  projective nonsingular  surface   with
$\kappa$-\emph{dim}$(Y) \ge 0$, and
 let $C\subset Y$ be a general   nonsingular connected curve of   genus $g\ge 3$
 moving in a positive dimensional linear system.
 Then
 \[
 g \le 6+5 p_g(Y)+ {1 \over 2}h^0(K_Y-C)
\]
  In particular:
\[
g \le \begin{cases}6,&\text{if  $p_g(Y)=0$} \\ 11,&\text{if
$p_g(Y)=1$}
 \end{cases}
\]
 \end{theorem}

 This theorem implies in particular the well known bound $g \le 11$ for a
 general nonsingular curve of genus $g$ on a K3--surface
 \cite{aM72}, and $g \le 6$ for a general nonsingular curve
 of genus $g$ on an Enriques surface.

 The  following result shows that   general curves on \emph{irregular}  surfaces
 of positive geometric genus cannot move in a positive dimensional linear
 system.

 \begin{theorem}\label{T:intr4}
Let $Y$ be a projective nonsingular surface with  $p_g>0$ and
$q>0$, and  let $C\subset Y$ be a nonsingular curve of genus $g
\ge 3$ such that  $\dim(|C|) \ge 1$.  Then no pencil $\Lambda
\subset |C|$ containing $C$ as a member defines a free fibration.
\end{theorem}

Theorem \ref{T:intr3}  does not say much  when the geometric
genus $p_g(Y)$ is large, and therefore it is ineffective for most
surfaces of general type. Our main result in this case is the
following.

 \begin{theorem}\label{T:intr1}
Let  $Y$ be a   surface of general type and let  $Z$ be the
minimal model of $Y$. Assume that $K_Z^2 \ge 3\chi(\O_Z) - 10$ and
that $C \subset Y$  is a general nonsingular  connected curve of
genus $g \ge 3$.   If     one of the following holds:
\begin{itemize}
\item[(a)] $\dim(|C|) \ge 2$,

\item[(b)]  $\dim(|C|)=1$, $h^0(K_Y-C)=0$ and  $C^2\ge {(g-1)\over
2}$,

\item[(c)] $\dim(|C|)=1$, and $h^1(C,\O_C(2C))=0$.
\end{itemize}
then
 $g \le 19$.
 \end{theorem}

Recall that  $K_Z^2 \ge 3\chi(\O_Z) - 10$ is the \emph{Castelnuovo
inequality} and if it is violated then  $Z$ is a double cover of a
ruled surface $R$. In this case, if $Y$ contains a general curve
such that $\dim(|C|) \ge 1$  then $C$ is mapped birationally into
$R$ and it follows that $R$ is rational. I have not been able to
obtain a meaningful bound for $g$ in this case.

The restrictions (b),(c) in the statement of the theorem might a
priori be sharp, and not merely  due to the method of proof: as
the genus increases, it seems difficult to  find a smooth curve
with general moduli that moves in a pencil with small
self-intersection on a surface.

The paper is divided into sections as follows.  In \S 1 we
introduce the terminology and main calculations concerning
fibrations. In \S \ref{S:ratfib} we develop the deformation theory
of fibrations.   Free fibrations are introduced in \S
\ref{S:free}. In the next \S \ref{S:freegen}   general curves are
studied and a characterization of the uniruledness of $\mgbar$ is
given.  The basic inequality  (\ref{E:ineq2})  is proved in \S
\ref{S:estim}. In \S  \ref{S:surf1} we prove Theorems
\ref{T:intr2}, \ref{T:intr3} and \ref{T:intr4}, while in \S
\ref{S:surfgt} we prove Theorem \ref{T:intr1}.

 \medskip

 \emph{Acknowledgements} - I  thank E. Ballico, L. Benzo, A. Bruno, F. Catanese,
  C. Ciliberto, A. Lopez,  M. Reid, M. Roth and A. Verra
  for useful conversations related to this work. Remarks and questions of the referee
	 contributed to improve the paper significantly. I am very thankful to him.

 \section{Fibrations}

 If
\[
\varphi:W \to U
\]
is a morphism of algebraic schemes, and $u\in U$, we denote by
$W(u)$ the scheme-theoretic fibre of $\varphi$ over $u$.

A \emph{$(-1)$-curve} (resp. a \emph{$(-2)$-curve}) on a
projective nonsingular surface $Z$ is a nonsingular connected
rational curve with self-intersection $-1$ (resp. $-2$). We will
denote by $K_Z$ a canonical divisor on   $Z$. By a \emph{minimal
surface} we   mean, as customary,  a projective nonsingular
connected surface without $(-1)$-curves.

By a \emph{fibration} we mean a surjective morphism
$$
f:  X \to S
$$
with connected fibres from a projective nonsingular surface to a
projective nonsingular connected curve. We will always denote by
\begin{itemize}

\item[]  $g =$ \quad the genus of the general fibre. \emph{We will
always assume} $g \ge 2$.

\item[]  $b =$ \quad  the genus of $S$.

\end{itemize}

A fibration  is called:

\begin{itemize}

\item[-] \emph{relatively minimal} if there are no $(-1)$-curves
contained in any of its fibres.

\item[-]  \emph{semistable} if it is relatively minimal and  every
fibre   has at most nodes  as singularities.

\item[-] \emph{stable} if it is relatively minimal and  every
fibre  is a stable curve (i.e. at most ordinary nodal singularities and with finitely many automorphisms).

\item[-] \emph{isotrivial} if all of its nonsingular fibres are
mutually isomorphic; equivalently, if two general nonsingular
fibres of $f$ are mutually isomorphic (the equivalence of the two
formulations follows from the separatedness of $\mgbar$).

\end{itemize}

Sheaves of differentials will be denoted with the symbol $\Omega^1$
and dualizing sheaves with the symbol $\omega$. Let  $f: X \to S$
be any fibration.
   Since both $X$ and $S$ are nonsingular we have:
$$
\wedge^2\Omega^1_X = \omega_X,  \quad \Omega^1_S = \omega_S
$$
Moreover $f$ is a relative complete intersection morphism and it
follows that
\begin{equation}\label{E:setup0}
\omega_{X/S} = \omega_X\otimes f^*\omega_S^{-1}
\end{equation}
(\cite{sK80}, Corollary (24)). We have an   exact sequence
\begin{equation}\label{E:setup3}
0 \to f^*\omega_S \to \Omega^1_X \to \Omega^1_{X/S} \to 0
\end{equation}
(which is exact   on the left because the first homomorphism is
injective on a dense open set and $f^*\omega_S$ is locally free).
If we dualize the sequence (\ref{E:setup3})   we obtain the exact
sequence:
\begin{equation}\label{E:setup6}
0 \to  T_{X/S}   \to T_X \to f^* T_S \to N\to 0
\end{equation}
where we have denoted
$$
T_{X/S} := Hom_{\O_X}(\Omega^1_{X/S},\O_X)
$$
and
\begin{equation}\label{E:setup2.1}
   N := Ext_{\O_X}^1(\Omega^1_{X/S},\O_X)
\end{equation}
The sequence (\ref{E:setup6}) shows, in particular, that $N$ is
the normal sheaf of $f$, and also the first relative cotangent
sheaf $T^1_{X/S}$ of  $f$  (see \cite{eS06}).   In particular,
\emph{$N$ is supported on the set of singular points of the fibres
of $f$. }

Moreover, as remarked in \cite{fS92}, p. 408, \emph{$T_{X/S}$ is
an invertible sheaf}  because it is a second syzygy of the
$\O_X$-module  $N$.

\begin{lemma}\label{L:setup1}
If all the singular  fibres of $f$ are reduced then
\begin{equation}\label{E:setup5}
T_{X/S} \cong \omega_{X/S}^{-1}
\end{equation}
\end{lemma}

\proof Since  the locus where $f$ is not smooth is finite,
$\Omega^1_{X/S}$ and $\omega_{X/S}$ are isomorphic in codimension
one. We thus have
\[
 c_1(\Omega^1_{X/S}) = c_1(\omega_{X/S})
 \]
 and the conclusion follows because $T_{X/S}$ is invertible.  \qed

 \begin{remark}\rm
 If $f$ has some non-reduced  singular fibre then   $T_{X/S}$ and $\omega_{X/S}^{-1}$ are
  not isomorphic. Precisely we have:
 \[
 T_{X/S}  = \omega_{X/S}^{-1}(\sum (\nu_i-1)E_i)
 \]
 where $\{\nu_iE_i\}$ is the set of all components of the singular fibres of $f$, and where $\nu_i$
 denote their corresponding multiplicities.  This formula is due to Serrano (\cite{fS92}, Lemma 1.1).
 \end{remark}

 The next proposition  will not be applied in the sequel:

 \begin{prop}\label{P:setup1}
If $f$ is any fibration we have:
$$
\chi(\omega_{X/S}^{-1}) = \chi(\O_X) +K_X^2 - 6(b-1)(g-1)
$$
If  the singular fibres of  $f$ are reduced  then
\[
\chi(T_{X/S})=\chi(\O_X) + K_X^2 - 6(b-1)(g-1)
\]
\end{prop}

\proof By Riemann--Roch:
$$
\chi(\omega_{X/S}^{-1}) = \chi(\O_X)+ {1\over
2}[\omega_{X/S}\cdot\omega_{X/S}+\omega_X\cdot\omega_{X/S}]
$$
(replacing $\omega_X \otimes f^*\omega_S^{-1}$ for $\omega_{X/S}$)
$$
\begin{array}{l}
= \chi(\O_X)+{1\over 2}[K_X^2 - 2(2b-2)(2g-2) +K_X^2 -(2b-2)(2g-2)]  \\\\
= \chi(\O_X)+ K_X^2 - 6(b-1)(g-1)
\end{array}
$$
The last formula is a consequence of Lemma \ref{L:setup1}. \qed

 The following is a classical result of Arakelov in the semistable case, and it is due to
 Serrano in the general case  (recall that we are assuming $g \ge 2$):

\begin{theorem} \label{T:setup1}
If $f$ is a non-isotrivial  fibration   then
\begin{equation}\label{E:setup4}
h^0(X,T_{X/S})  = h^0(X,T_X) = 0
\end{equation}
If moreover $f$ is relatively minimal then we also have:
\begin{equation}\label{E:setup4.1}
h^1(X,T_{X/S}) = 0
\end{equation}
\end{theorem}

 \proof
 $f$ is non-isotrivial if and only if $f_*T_X=0$ (\cite{fS96}, Lemma 3.2). Since
 \[
 f_*T_{X/S} \subset f_*T_X
 \]
 we also have  $f_*T_{X/S}=0$ if $f$ is non-isotrivial. Thus (\ref{E:setup4})
 is a consequence of the Leray spectral sequence.
 For (\ref{E:setup4.1}) see  \cite{sA71} or \cite{lS81} in the semistable case,
 and  \cite{fS92}, Corollary 3.6, in the general case. \qed

     Denoting by  $Ext^1_f$  the first derived functor of $f_*Hom$,
  we are interested in the sheaf  $Ext^1_f(\Omega^1_{X/S},\O_X)$ because its cohomology controls
  the deformation theory of $f$ (see Lemma \ref{L:free0} below).  This sheaf is not locally free in general, but it decomposes 
	as follows:
	\begin{equation}\label{E:extdec}
	Ext^1_f(\Omega^1_{X/S},\O_X)  = \mathcal{E} \oplus \mathcal{T}
\end{equation}
where $\mathcal{E}$ is locally free and $\mathcal{T}$ is a torsion sheaf. By the \emph{rank}  of 
$Ext^1_f(\Omega^1_{X/S},\O_X)$ 
we will mean the rank of $\mathcal{E}$.

  \begin{lemma}\label{L:setup2}
  For any non-isotrivial fibration $f:X \to S$ 
   the   sheaf $Ext^1_f(\Omega^1_{X/S},\O_X)$ has rank $3g-3$. Moreover there is
    an exact sequence of sheaves on $S$:
  \begin{equation}\label{E:dpsi2}
\begin{array}{ccccc}
0 &\to R^1f_*T_{X/S}& \buildrel c_{10}\over\longrightarrow
 Ext^1_f(\Omega^1_{X/S},\O_X) \buildrel c_{01}\over\longrightarrow
& f_* Ext^1_X(\Omega^1_{X/S},\O_X)&\to 0   \\
& &&\Vert \\
& &&f_*N
\end{array}
\end{equation}
If   $Ext^1_f(\Omega^1_{X/S},\O_X)$ is locally free  and all the fibres of $f$ are reduced 
then there is an isomorphism
  \begin{equation}\label{E:dual1}
  Ext^1_f(\Omega^1_{X/S},\O_X) \cong Hom(f_*(\Omega^1_{X/S}\otimes\omega_{X/S}),\O_S)
  \end{equation}
 \end{lemma}

 \proof
 If $p\in S$ is such that $X(p)$   is smooth then  $Ext^1(\Omega^1_{X(p)},\O_{X(p)})= H^1(T_{X(p)})$ has dimension $3g-3$. 
Then, if $U \subset S$ is the open set over which $f$ is smooth,  
$Ext^1_f(\Omega^1_{X/S},\O_X)_{|U}$
is locally free of rank $3g-3$.  

(\ref{E:dpsi2})  is the sequence associated to the local-to-global spectral sequence for $Ext_f$.

Assume that $Ext^1_f(\Omega^1_{X/S},\O_X)$ is locally free and that all the fibres of $f$ are reduced. Then, since 
\[
Ext^0_f(\Omega^1_{X/S},\O_X) = f_*Hom(\Omega^1_{X/S},\O_X) = f_*T_{X/S}=0
\]
(see the proof of Theorem \ref{T:setup1}) is   locally free as well, and the higher $Ext^i_f$'s vanish, both sheaves 
commute with base change (\cite{hL83}, Th. 1.4). The reducedness of the fibres implies that $\Omega^1_{X/S}$ is torsion free,
hence flat over $S$. Therefore (\ref{E:dual1})   follows from  relative duality  for the fibration $f$ (\cite{sK80}, Corollary (24)). \qed

\begin{remark}\label{R:locfree}\rm
The torsion sheaf $\mathcal{T}$ in (\ref{E:extdec}) is non-zero in general, even if all the fibres of the non-isotrivial fibration $f$ are reduced. Assume for example that for some $p \in S$ there is a $(-2)$-curve in the fibre $X(p)$. Then $\mathrm{Hom}(\Omega^1_{X(p)},\O_{X(p)}) \ne 0$ because $X(p)$ has non-trivial infinitesimal automorphisms. Since
$f_* Hom(\Omega^1_{X/S},\O_X)=f_*T_{X/S}=0$ this sheaf does not commute with base change and therefore $Ext^1_f(\Omega^1_{X/S},\O_X)$ 
cannot be locally free (\cite{hL83}, Th. 1.4).
A fortiori there is torsion if the fibration is not relatively minimal.
If the fibration $f$ is stable then this phenomenon does not occur and $\mathcal{T}=0$.

On the other hand, if $f$ is relatively minimal and non-isotrivial then $R^1f_*T_{X/S}$ is locally free
even if there are $(-2)$-curves in the fibres or non-reduced fibres.
This follows from the fact that  $H^0(R^1f_*T_{X/S})=0$ 
(because $H^1(X,T_{X/S})=0$  by Theorem \ref{T:setup1}). 
\end{remark}

\begin{prop}\label{P:def1}
If the fibration $f$ is non-isotrivial  then we have:
\begin{equation}\label{E:dpsi1}
\chi(Ext^1_f(\Omega^1_{X/S},\O_X)) = 11\chi(\O_X) - 2 K_X^2 +
2(b-1)(g-1)
\end{equation}
\end{prop}

  \proof
 Since the fibres of $f$ are 1-dimensional we have
 \[
 R^2f_*T_{X/S} = 0
 \]
 Moreover $f_*Ext^2(\Omega^1_{X/S},\O_X)=0$ because $Ext^2(\Omega^1_{X/S},\O_X)=0$
  by the exact sequence (\ref{E:setup3}).
 Therefore, using the local-to-global spectral sequence for $Ext_f$ we deduce that
 \[
R^1 f_*N = Ext^2_f(\Omega^1_{X/S},\O_X) = 0
\]
 where the last equality is true because the fibres of $f$ are 1-dimensional.  This gives:
 \[
 \chi(f_*N)=\chi(N)
 \]
 Moreover, since $f$ is non-isotrivial, from (\ref{E:setup4}) and the Leray spectral sequence we get
 \[
  \chi(R^1f_*T_{X/S}) = - \chi(T_{X/S})
  \]
   We now use the exact sequence (\ref{E:dpsi2}) and we deduce that
  \[
  \begin{array}{llll}
  \chi(Ext^1_f(\Omega^1_{X/S},\O_X)) &= \chi(f_*N)+\chi(R^1f_*T_{X/S}) \\
  &=\chi(N)-\chi(T_{X/S})  \\
  &= \chi(f^*T_S) - \chi(T_X)& (\hbox{by (\ref{E:setup6})}) \\
  \end{array}
  \]
   Using Riemann-Roch one computes that:
   \[
   \begin{array}{ll}
   \chi(f^*(T_S)) &= \chi(\O_X)+ 2(b-1)(g-1)  \\
   \chi(T_X) & = 2K_X^2 - 10\chi(\O_X)
   \end{array}
   \]
   and by substitution one gets (\ref{E:dpsi1}).   \qed

 \begin{remark}\rm
   Assume that the fibration $f$ is   stable.
   Then with $f$ there is associated a  modular morphism
$$
\psi_f: S \longrightarrow \overline{\M}_g
$$
to the moduli stack of  stable curves of genus $g$.  We know
(\cite{HM82}, p. 49) that:
$$
f_*(\Omega^1_{X/S}\otimes\omega_{X/S}) = \psi_f^*\Omega^1_
{\overline{\M}_g}
$$
From (\ref{E:dual1})  it follows that:
\begin{equation}\label{E:dpsi4.1}
Ext^1_f(\Omega^1_{X/S},\O_X) =  \psi_f^*T_{\overline{\M}_g}
\end{equation}
Therefore the identity (\ref{E:dpsi1}) can be recovered by
applying
  the Riemann-Roch theorem to the vector bundle
  $\psi_f^*T_{\overline{\M}_g}$ as follows.
Recall that
\begin{equation}\label{E:dpsi4}
 c_1(\psi_f^*K_{\overline{\M}_g})= 13\lambda(f)-2\delta(f) = -5(b-1)(g-1) - [11\chi(\O_X) - 2K_X^2]
 \end{equation}
where $\lambda(f) = \chi(\O_X)-(b-1)(g-1)$ and $\delta(f)=
12\chi(\O_X)-K_X^2-4(b-1)(g-1)$.  Then:
\[
\begin{array}{ll}
\chi( \psi_f^*T_{\overline{\M}_g})&=
 -c_1(\psi_f^*K_{\overline{\M}_g})+(1-b)(3g-3)\\
&=11\chi(\O_X) - 2 K_X^2 + 2(b-1)(g-1) \end{array}
\]
 that is
(\ref{E:dpsi1}).

In case the fibration $f$ is not   stable we still have a
non-empty open set $U \subset S$ above which all fibres of $f$ are
  stable. Therefore we have an induced morphism $U
\longrightarrow \stgbar$; since $S$ is a nonsingular curve this
morphism extends to a morphism:
\[
\overline{\psi}_f: S \longrightarrow \mgbar
\]
with values in the moduli \emph{space.}    Now we no longer have
the interpretation (\ref{E:dpsi4.1}), even replacing $\stgbar$ by
$\mgbar$. This will be the reason for some lengthening  in the
following discussion of the deformation theory of $f$.
\end{remark}

\section{Deformation theory}\label{S:ratfib}

\begin{lemma}\label{L:free0}
Let  $f:X \to S$ be a  non-isotrivial  fibration. Then there is a
natural isomorphism
\[
\mu: \mathrm{Ext}^1_X(\Omega^1_{X/S},\O_X) \to H^0(S,
Ext^1_f(\Omega^1_{X/S},\O_X))
\]
and both spaces are naturally identified with the tangent space of
$\mathrm{Def}_f$,
 the functor of Artin rings of deformations of $f$ leaving the target fixed (see \cite{eS06}, p. 164).
Moreover
\[
  H^1(S, Ext^1_f(\Omega^1_{X/S},\O_X))
\]
is an obstruction space for $\mathrm{Def}_f$.
\end{lemma}

\proof There is an exact sequence:
\[
0 \to H^1(S,f_*T_{X/S}) \to \mathrm{Ext}^1_X(\Omega^1_{X/S},\O_X)
\buildrel\mu\over\longrightarrow
 H^0(S, Ext^1_f(\Omega^1_{X/S},\O_X)) \to 0
 \]
 (\cite{hL83}, p. 105). Since $f_*T_{X/S}=0$ (see the proof of Theorem \ref{T:setup1})  $\mu$ is an isomorphism.
First order deformations of $f$ are in 1-1 correspondence with
isomorphism classes of extensions
\[
\zeta: 0 \to \O_X \buildrel j\over\longrightarrow \O_{\X} \to \O_X
\to 0
\]
such that  $\O_{\X}$ is a   sheaf of flat
$\O_S[\epsilon]$-algebras.   In \cite{eS06}, Theorem 1.1.10,
 it is proved that if we only assume that $\zeta$  is an extension of $\O_S$-algebras then it has  already
 a structure of deformation by sending $\epsilon \mapsto j(1)$.
  There is a natural 1-1 correspondence between isomorphism classes of extensions $\zeta$
  as above and  $\mathrm{Ext}^1_X(\Omega^1_{X/S},\O_X)$ (loc. cit., Theorem 1.1.10).

By comparing the cohomology sequence of (\ref{E:dpsi2}) with the
exact sequence (d) of Lemma 3.4.7 in loc. cit., we see that
$H^1(S,Ext^1_f(\Omega^1_{X/S},\O_X))$ is an obstruction space for
$\mathrm{Def}_f$. \qed

Assume given a non-isotrivial  fibration $f:X\longrightarrow S$.  
Given $p \in S$   we have a morphism of functors:
\[
\psi_p:  \mathrm{Def}_f \longrightarrow  \mathrm{Def}_C
\]
where $C:=X(p)$, defined as follows. Given a local artinian $\mathbb{C}$-algebra $A$ and 
an element   $ \eta_A\in\mathrm{Def}_f(A)$:
\[
\xymatrix{
X \ar@{^(->}[r]\ar[d]^f &X_A \ar[d]^{f_A} \\
    S \ar@{^(->}[r]_-{\mathrm{id}_S\times \{v\}} & S\times \mathrm{Spec}(A)}  
		\]
we define $\psi_p(\eta_A)\in \mathrm{Def}_C$ to be the left square of 
 \[
    \xymatrix{
 C\ar[r]\ar[d]&C_A \ar[d]\ar@{^(->}[r]&X_A \ar[d]^{f_A} \\
\mathrm{Spec}(\mathbb{C})\ar[r]^-{\{p\}}&\mathrm{Spec}(A)\ar[r]^-{\{p\}\times 1_A} & S\times\mathrm{Spec}(A) }
\]

\begin{lemma}\label{L:basech}
Let $f:X\longrightarrow S$ be a non-isotrivial  fibration and $p
\in S$ such that the fibre $X(p)$ is reduced. Then:
\begin{description}
\item[(i)] There is an exact sequence of sheaves on $S$:
\begin{footnotesize}
\begin{equation}\label{E:basech}
    \xymatrix{
    0\ar[r]&Ext^1_f(\Omega^1_{X/S},\O_X)(-p)\ar[r]&Ext^1_f(\Omega^1_{X/S},\O_X)\ar[r]&
    \mathrm{Ext}^1_{X(p)}(\Omega^1_{X(p)},\O_{X(p)})_p\ar[r]&0}
\end{equation}
\end{footnotesize}
where $\mathrm{Ext}^1_{X(p)}(\Omega^1_{X(p)},\O_{X(p)})_p$ is the
skyscraper sheaf supported on $p$ with fibre
$\mathrm{Ext}^1_{X(p)}(\Omega^1_{X(p)},\O_{X(p)})$.

\item[(ii)] The homomorphism $\kappa_f:T_S \longrightarrow
Ext^1_f(\Omega_{X/S},\O_X)$ defined by the exact sequence
(\ref{E:setup3}) is injective and the induced map:
\[
\kappa_f(p): T_pS
\longrightarrow\mathrm{Ext}^1_{X(p)}(\Omega^1_{X(p)},\O_{X(p)})
\]
  coincides with the Kodaira-Spencer map of the family $f$ at
$p$.

\item[(iii)] Let
\[
\xymatrix{r_p:H^0(S,Ext^1_f(\Omega^1_{X/S},\O_X))\ar[r]&
\mathrm{Ext}^1_{X(p)}(\Omega^1_{X(p)},\O_{X(p)})}
\]
 be  the map induced by (\ref{E:basech}). Then $r_p$ is the differential of $\psi_p$.

  \item[(iv)]  If $S=\P^1$ then
  $\mathrm{Im}(\kappa_f(p)) \subset \mathrm{Im}(r_p)$.
 \end{description}
\end{lemma}

\proof (i)
 Consider the base
change map:
\[ \xymatrix{\tau^1(p):Ext^1_f(\Omega^1_{X/S},\O_X)_p\otimes\k(p)
  \ar[r]&\mathrm{Ext}^1_{X(p)}(\Omega^1_{X(p)},\O_{X(p)})}\]
  Then the base change theorem for the relative Ext sheaves and
  the fact that $Ext^2_f(\Omega^1_{X/S},\O_X)=0$ imply that
  $\tau^1(p)$ is an isomorphism for all $p \in S$ (\cite{hL83}, Th. 1.4).  Therefore
  (\ref{E:basech}) is the exact sequence obtained by tensoring the
  sequence
  \[
  \xymatrix{0\ar[r]&\O_S(-p) \ar[r]&\O_S \ar[r] & \O_p \ar[r] &0}\]
by $Ext^1_f(\Omega^1_{X/S},\O_X)$.

(ii) The exact sequence (\ref{E:setup3}) induces an exact sequence
on $S$:
 \[
0 \to f_*T_X \to T_S  \buildrel\kappa_f\over\longrightarrow
Ext^1_f(\Omega^1_{X/S},\O_X)
 \]
 Since $f$ is non-isotrivial we have $f_*T_X=0$
  by Lemma 3.2 of \cite{fS96}, and this proves the first assertion. $\kappa_f(p)$ is the composition
\[
 \xymatrix{T_pS \ar[r]&
 Ext^1_f(\Omega^1_{X/S},\O_X)_p\otimes\k(p)\ar[r]^-{\tau^1(p)}&
 \mathrm{Ext}^1_{X(p)}(\Omega^1_{X(p)},\O_{X(p)})}
\]
 of the natural restriction of $\kappa_f$ with the base change map $\tau^1(p)$.
 Since $X(p)$ is reduced   $\mathrm{Ext}^1_{X(p)}(\Omega^1_{X(p)},\O_{X(p)})$ is the  vector space
 of first order deformations of $X(p)$ and $\kappa_f(p)$ is the Kodaira-Spencer map of $f$ at $p$
 essentially by definition (see \cite{eS06} , Remark 2.4.4).

(iii) is tautological.

(iv) If $S=\P^1$ then every $\theta\in T_p\P^1$ comes from a
global vector field $\Theta \in H^0(T_{\P^1})$. Then the assertion follows from the commutativity of the diagram:
 \[
\xymatrix{
H^0(T_{\P^1})\ar[r]\ar[d] & H^0(Ext^1_f(\Omega^1_{X/S},\O_X)) \ar[d]^-{r_p} \\
T_p\P^1 \ar[r]^-{\kappa_f(p)} & \mathrm{Ext}^1_{X(p)}(\Omega^1_{X(p)},\O_{X(p)})}
\]
 \qed

Every non-isotrivial   fibration
$f$ has a semiuniversal formal deformation  $(R, \{\eta_n\})$, where 
$\eta_n \in \mathrm{Def}_f(R/\mathbf{m}^{n+1})$, $n \ge 1$, and $\mathbf{m}\subset R$
is the maximal ideal of the local complete $\mathbb{C}$-algebra $R$ (\cite{eS06}, Th. 3.4.8).
The next result states the existence of an \emph{algebraic} semiuniversal deformation in the relatively minimal case.

\begin{theorem}\label{T:semiuniv}
Let $f:X\longrightarrow S$ be a non-isotrivial  relatively minimal fibration. Then $f$ has 
a    semiuniversal algebraic deformation, i.e. a deformation:
\[
    \xymatrix{
    X \ar@{^(->}[r]\ar[d]^f &\X \ar[d]^F \\
    S \ar@{^(->}[r]_-{\mathrm{id}\times \{v\}} & S\times V }
\]
 parametrized by a pointed algebraic scheme $(V,v)$ which is formally semiuniversal at $v$.
\end{theorem}

\proof The functor $\mathrm{Def}_{(X,\omega_{X/S})}$ of deformations of the  pair $(X,\omega_{X/S})$ 
has a formal semiuniversal deformation $(\O, \{(X_n,L_n)\})$, where $\O$ is a complete local $\mathbb{C}$-algebra with maximal ideal
$\mathbf{I}$ and $(X_n,L_n)$ is a deformation of $(X,\omega_{X/S})$ over $\O/\mathbf{I}^{n+1}$ (\cite{eS06}, Th. 3.3.11). 
Let $Y$ be the surface obtained by contracting all the $(-2)$-curves of $X$ contained in the fibres of $f$. From Theorem 2' of 
\cite{lS81} and the Nakai-Moisezon-Kleiman criterion it follows that $\omega_{X/S}$ is the pullback of an ample invertible sheaf on $Y$. 
Therefore a positive power 
 $\omega_{X/S}^{\otimes k}$  is globally generated and maps $X$ birationally to a projective surface whose only singularities are 
rational double points. Arguing as in \cite{mA74}, Example 5.5 one deduces the existence of 
a  semiuniversal algebraic deformation of the pair $(X,\omega_{X/S})$, i.e. of a pair 
consisting of   a    deformation of $X$:
 \[
 \xymatrix{
 X \ar@{^(->}[r]\ar[d] &\mathcal{Z} \ar[d]^\zeta \\
 \mathrm{Spec}(\mathbb{C})\ar@{^(->}[r]_-{ \{w\}} &W}
 \]
 parametrized by a pointed algebraic scheme $(W,w)$ and of an invertible sheaf $\L$ on $\mathcal{Z}$ such that the pair $(\zeta, \L)$ is a deformation of 
$(X,\omega_{X/S})$ which is formally semiuniversal at $w$.  Now   let   $(V,v)$   be the relative 
local Hom-scheme $Hom(\mathcal{Z}/W,S\times W/W)$ around the point $[f]$ representing
 $f: X=\mathcal{Z}(w)\longrightarrow S\times \{w\}$ (see
 \cite{jK99}, Ch. I.1).   Then it is straightforward to check that the family of deformations of $f$ parametrized by $(V,v)$ has the required properties. \qed

\section{Free fibrations}\label{S:free}

From now on we will only consider fibrations parametrized by
$\P^1$. All examples  of  such fibrations are obtained as follows.

Let   $Y$ be a projective nonsingular surface,
 and let $C \subset Y$ be a projective nonsingular connected curve of genus $g$ such that
\begin{equation}\label{E:posdim}
\dim(|C|)   \ge 1
\end{equation}
Consider a linear pencil $\Lambda$ contained in $|C|$ whose
general member is nonsingular and let   $\sigma: X\to Y$ be the
blow-up at its  base points (including the infinitely near ones).
We obtain a      fibration
$$
f: X \to \P^1
$$
by choosing an isomorphism $\Lambda^\vee\cong \P^1$ and taking the
composition:
\[
\xymatrix{X\ar[r]^-\sigma & Y\ar@{-->}[r]&
\Lambda^\vee\cong \P^1}
\]
 We will call $f$  \emph{the fibration
defined by the pencil $\Lambda$}.

In the decomposition (\ref{E:extdec}) we have  
\begin{equation}\label{E:freea}
\mathcal{E} = \bigoplus_{i=1}^{3g-3}
\O_{\P^1}(a_i)
\end{equation}
for some integers $a_i$.

\begin{definition}\label{D:free1}
Let $f: X \to \P^1$ be a  fibration.    We call $f$ \emph{free}  if it is
non-isotrivial  and $Ext^1_f(\Omega^1_{X/\P^1},\O_X)$ is globally
generated, i.e. if in (\ref{E:freea}) we have $a_i \ge 0$ for all $i$.  
\end{definition}

 \begin{prop}\label{P:def2}
 Suppose that $f:X \to \P^1$ is a  free fibration. Then:
 \begin{description}
    \item[(i)] The functor  $\mathrm{Def}_f$ is smooth of
    dimension $\ge 3g-1$.
    \item[(ii)]
 \begin{equation}\label{E:free1}
 11\chi(\O_X) - 2 K_X^2 \ge  5(g-1)+2
 \end{equation}
 \end{description}
 
   \end{prop}

\proof
  (i) Since $a_i \ge 0$ for all $i$ in (\ref{E:freea}),
  we    have  $h^1(\P^1,Ext^1_f(\Omega^1_{X/\P^1},\O_X)) =0$.
  Therefore $\mathrm{Def}_f$ is smooth, by Lemma  \ref{L:free0}.
    Since $f$ is non-isotrivial   the homomorphism
  \[
  \kappa_f:T_{\P^1} \to Ext^1_f(\Omega^1_{X/\P^1},\O_X)
  \]
 induced by   the exact sequence (\ref{E:setup3})  is injective, by Lemma \ref{L:basech}.
  It follows that in (\ref{E:freea})  we have   $a_i \ge 2$ for some $i$
 and therefore $\mathrm{Def}_f$ has dimension
 \begin{equation}\label{E:free2}
  h^0(\P^1,Ext^1_f(\Omega^1_{X/\P^1},\O_X)) \ge 3g-1
 \end{equation}

 (ii)
    Recalling  Proposition  \ref{P:def1} we obtain:
  $$
  \begin{array}{llll}
 11\chi(\O_X) - 2 K_X^2 - 2(g-1)&=\chi(Ext^1_f(\Omega^1_{X/\P^1},\O_X))  \\ \\
 &= h^0(\P^1, Ext^1_f(\Omega^1_{X/\P^1},\O_X)) \\ \\
 &\ge  3g-1
 \end{array}
 $$
 and this is equivalent to  (\ref{E:free1}). \qed

\begin{remark}\label{R:free1}\rm
If $f$ is   stable then recalling (\ref{E:dpsi4}) we see that
(\ref{E:free1}) is equivalent to the condition
\[
c_1(\psi_f^*K_{\mgbar}) \le -2
\]
 \end{remark}

Consider now a non-isotrivial fibration $f:X \longrightarrow \P^1$
(not necessarily   stable anymore) and the associated exact
sequence:
\begin{equation}\label{E:reldiff1}
    \xymatrix{
    0 \ar[r]& f^*\omega_{\P^1} \ar[r]& \Omega^1_X \ar[r]& \Omega^1_{X/\P^1} \ar[r]& 0
    }
\end{equation}
Taking $f_*Hom(-,\O_X)$ of it and recalling that $f_*T_X=0$
(\cite{fS96}, Lemma 3.2) we obtain the following 4-term exact
sequence of sheaves on $\P^1$:
\begin{equation}\label{E:reldiff2}
\xymatrix{ 0 \ar[r]&
T_{\P^1}\ar[r]&Ext^1_f(\Omega^1_{X/\P^1},\O_X)\ar[r]& R^1f_*T_X
\ar[r]& T_{\P^1}\otimes R^1f_*\O_X \ar[r] &0}
\end{equation}
We thus obtain a long exact cohomology sequence as follows:
\begin{small}\begin{equation}\label{E:reldiff3}
\begin{array}{ccccc}
&&H^1(T_X) \\
&&|| \\
 0 \longrightarrow
H^0(T_{\P^1})\longrightarrow&H^0(Ext^1_f(\Omega^1_{X/\P^1},\O_X))\longrightarrow&
H^0(R^1f_*T_X)\to \\ \\
  \to H^0(T_{\P^1}\otimes
R^1f_*\O_X)\longrightarrow&H^1(Ext^1_f(\Omega^1_{X/\P^1},\O_X))\longrightarrow&
H^1(R^1f_*T_X) \longrightarrow&H^1(T_{\P^1}\otimes
R^1f_*\O_X)\longrightarrow 0 \\
&&|| \\
&&H^2(T_X)
\end{array}
\end{equation}
\end{small}

\begin{lemma}\label{L:blowup}
Let $g:Y \longrightarrow \P^1$ be a non-isotrivial fibration,  $\pi: X \longrightarrow Y$  the blow-up
of a point $y \in Y$ and $f = g \cdot \pi:X \longrightarrow \P^1$.  Then   $f$ is free if and only if $g$ is free.
\end{lemma}

\proof Let $\pi^{-1}(y)=E$ be the exceptional curve. Then there is an exact sequence (\cite{eS06}, p. 172)
\[
\xymatrix{
0 \ar[r] & T_X \ar[r] & \pi^* T_Y \ar[r] & \O_E(1) \ar[r] & 0}
\]
which implies that we have a surjection $\epsilon: R^1f_*T_X \longrightarrow R^1g_*T_Y$. 
 Now compare the exact sequence (\ref{E:reldiff2}) with the analogous one for $g$:
\[
\xymatrix{ 0 \ar[r]&
T_{\P^1}\ar@{=}[d]\ar[r]&Ext^1_f(\Omega^1_{X/\P^1},\O_X)\ar@{-->}[d]\ar[r]& R^1f_*T_X
\ar[r]\ar[d]^-\epsilon& T_{\P^1}\otimes R^1f_*\O_X \ar[r]\ar@{=}[d] &0 \\
  0 \ar[r]&
T_{\P^1}\ar[r]&Ext^1_g(\Omega^1_{Y/\P^1},\O_Y)\ar[r]& R^1g_*T_Y
\ar[r]& T_{\P^1}\otimes R^1g_*\O_Y \ar[r] &0}
\]
This diagram implies the existence of the dotted arrow, which is necessarily surjective. 
Then $Ext^1_f(\Omega^1_{X/\P^1},\O_X)$ is globally generated if and only if $Ext^1_g(\Omega^1_{Y/\P^1},\O_Y)$ is.
 \qed

\section{Free fibrations and general curves}\label{S:freegen}

  We need to have  geometrical criteria
  to verify that a given fibration is free. For this purpose we introduce some
   definition. We start from the following well known one:

\begin{definition}\label{D:genmod}
A connected projective nonsingular curve $C$ of genus $g$
\emph{has general moduli,}  or \emph{is a general curve of genus
$g$}, if it is a general fibre in a smooth projective family
$\mathcal{C} \longrightarrow V$ of curves of genus $g$,
parametrized by a nonsingular connected algebraic scheme $V$, and
such that the  morphism $V\longrightarrow\mg$ induced by
functoriality is dominant. Such a family will be said to have
\emph{general moduli}.

\end{definition}

Whether a given family has general moduli can be checked in
practice by means of the Kodaira-Spencer map. The following are
some well known properties of a general curve that we will use:

\begin{prop}\label{P:gencur}
Let $C$ be a general   curve of genus $g\ge 3$.  Then:
\begin{description}
\item[(i)] $C$ has Clifford index
\[
{\rm Cliff}(C) = \bigg[{g-1 \over 2}\bigg]
\]
\item[(ii)] For any invertible sheaf $L$ on $C$, of degree $d\ge
0$ and with $h^0(L)=r+1 \ge 1$ we have
\[
\rho(g,r,d) \ge 0
\]
 where $\rho(g,r,d) := g - (r+1)(g-d+r)$ is the Brill-Noether number. In particular:
 \[
 \begin{array}{cc}
 d \ge {1 \over 2}g+1&\hbox{if $r\ge 1$} \\ \\
 d \ge {2 \over 3}g + 2 & \hbox{if $r\ge 2$}
 \end{array}
 \]
 \item[(iii)] $H^1(L^2)=0$ for  any invertible sheaf $L$ on $C$ such that $h^0(L)\ge 2$.

 \item[(iv)] $C$ does not possess irrational involutions, i.e. non-constant morphisms
   of degree $\ge 2$ onto a curve of positive genus.

   \item[(v)]  $C$ does not possess non-trivial automorphisms.

\end{description}

\end{prop}

\proof  (i)   follows from the main result of \cite{CM91} and from
Brill-Noether theory.
  (ii) is proved in \cite{GH80}. For (iii)    see
 \cite{AC81}, p. 22,  and references therein.  (iv) is classical and follows from Riemann's existence
 theorem and a dimension count.  For (v) see  \cite{hP69}.
\qed

The next definition we need is the following
  one:

\begin{definition}\label{D:gensurf}
Let $Y$ be a
  projective nonsingular   surface, and $C \subset Y$   a   projective nonsingular connected
curve   of genus $g \ge 3$. Let $r :=\mathrm{dim}(|C|)$ and denote
by $j:C\longrightarrow Y$ the embedding. We say that $C$ is a
\emph{general curve moving in an $r$-dimensional linear system on
$Y$}  if there is a family of deformations of $j$:
\begin{equation}\label{E:surf2}
\xymatrix{ C\ar[dr]\ar@{^{(}->}[dd]_j \ar@{^{(}->}[rr]&&\C \ar@{^{(}->}[dd]\ar[dr]   \\
&\mathrm{Spec}(\mathbb{C})\ar@{-}[r]&\ar[r]^{\{b\}} & B \\
Y\ar[ur]\ar@{^{(}->}[rr]&&\Y\ar[ur]_\beta}
\end{equation}
 parametrized by a pointed connected nonsingular algebraic
scheme $(B,b)$
 such that
 \begin{enumerate}
    \item the family of deformations of $C$ given by the upper  square has general moduli;
     \item   $\mathrm{dim}(|\C(u)|) \ge r$   on the surface $\Y(u)$ for all closed points $u \in B$.
 \end{enumerate}
 In the case $r=0$ we will just say that $C$ \emph{is a  curve
 with general moduli
 in} $Y$.
 \end{definition}

\begin{remark}\label{R:costable}\rm
In the case $r=0$ Definition \ref{D:gensurf} is equivalent to the
notion of \emph{costability} of $C$ in $Y$.  (see \cite{eH76} and
\cite{eS06}, Def. 3.4.22).
\end{remark}

\begin{example}\label{Ex:genus4}\rm
Let $C\subset \P^3$ be a nonsingular curve of type $(3,3)$ on a
nonsingular quadric $Q$.  It is easy to show that there is a
nonsingular quintic surface $Y$ containing $C$. On $Y$ we have
$C^2=0$ and dim$(|C|)=1$. Since $C$ is a canonical curve of genus 4
one can construct a family of deformations of the pair $(C,Y)$ so
that $C$ has general moduli. It follows that $C$ is a general
curve of genus 4 moving in a 1-dimensional linear system on $Y$.
This is a special case of a class of examples that can be
constructed  in a similar way. See also Example \ref{Ex:BV} and
\cite{aV11}, Example 2.3.
\end{example}

\begin{example}\label{Ex:spacecurves}\rm
Let $Y\subset \P^3$ be a nonsingular surface of degree $n \ge 6$,
$C \subset Y$ a  nonsingular curve of genus $g \ge 3$ such that
$\dim(|C|) \ge 1$. Then $L=\O_C(1)$ satisfies $h^0(L)\ge 3$ ($\ge
4$ if $C$ is non-degenerate) and $L^{n-4} =\omega_{Y|C}$ is
special. Therefore, by Proposition \ref{P:gencur}(iii), $C$ cannot
be a general curve. We therefore see that a general curve of genus
$g \ge 3$ cannot move in a positive dimensional linear system on a
nonsingular surface of degree $n \ge 6$ in $\P^3$.
\end{example}

   \begin{prop}\label{P:free1}
   Assume that
 $C$  is a general   nonsingular curve of genus $g \ge 3$  moving in a
 positive dimensional
linear system on a projective nonsingular  surface $Y$.  Then the
following conditions are equivalent:
\begin{description}
    \item[(i)] Each linear pencil $\Lambda\subset |C|$
    containing $C$ as a member defines  an isotrivial fibration.
    \item[(ii)] There is a linear pencil $\Lambda\subset |C|$
    containing $C$ as a member which defines  an isotrivial fibration.

    \item[(iii)] $Y$ is
birationally equivalent to $C\times \P^1$.
    \item[(iv)] $Y$ is a non-rational birationally ruled surface.

\end{description}

\end{prop}

   \proof $(i)\Longrightarrow (ii)$ is trivial.

$(ii)\Longrightarrow (iii)$.
 Let $f:X \to \P^1$ be the fibration
defined by the pencil $\Lambda$.
  Assume that $f$ is isotrivial. Then,  by the structure theorem for isotrivial fibrations  \cite{fS96},
  there is  a nonsingular curve $\Gamma$ and a finite group $G$
 acting on both  $C$ and  $\Gamma$ such that there is a birational isomorphism
 \[
 \xymatrix{X\ar@{-->}[r]&(C\times \Gamma)/G}
 \]
 and a commutative diagram:
 \[
 \xymatrix{X\ar[d]^-f \ar@{-->}[r]& (C\times \Gamma)/G\ar[d]
 \\ \P^1 \ar@{=}[r]& \Gamma/G}
 \]
 where the right vertical arrow is the projection.
   But since $C$ is general, it has no non-trivial automorphisms
   (Prop. \ref{P:gencur}(v)), and therefore $G$ acts trivially on $C$:
   thus $X$ is birational to $C \times (\Gamma/G)= C \times \P^1$.

$(iii)\Longrightarrow (iv)$ is obvious.

$(iv)\Longrightarrow (i)$.
   By hypothesis there is a  birational isomorphism
   $\xi: \xymatrix{Y \ar@{-->}[r]&  \Gamma\times \P^1}$ for some
    projective nonsingular curve $\Gamma$ of positive genus. Let $D$ be a general member of $\Lambda$.
    The composition
   \[
   \xymatrix{
   h: D \ \ar@{^{(}->}[r]& Y \ar@{-->}[r]^-\xi& \Gamma\times \P^1\ar[r] &\Gamma}
   \]
   where the last morphism is the projection, is non-constant.  Since $D$ is a general curve, it does not possess
    irrational involutions (Prop. \ref{P:gencur}(iv)), thus $h$ must be an
    isomorphism. Therefore all general  fibres of $f$ are mutually
    isomorphic, and this means that $f$ is isotrivial.  \qed

    The next proposition relates the notion of free fibration with
    Definition \ref{D:gensurf}.

\begin{prop}\label{C:free1}
 Assume that $C$  is
a general   nonsingular   curve  of genus $g \ge 3$ moving  in  a
positive-dimensional  linear system on a
  projective nonsingular  surface $Y$ which is not irrational ruled.  Then
 a general pencil $\Lambda \subset  |C|$ containing  $C$ as a member defines a free  fibration.
  \end{prop}

  \proof
  Let $f: X \longrightarrow \P^1$ be the fibration defined by $\Lambda$ and let $p \in \P^1$ be such that $C=X(p)$.
  By Proposition  \ref{P:free1},  $f$ is non-isotrivial.
By hypothesis there is a pointed nonsingular algebraic  scheme
$(B,b)$ and a commutative  diagram as follows:
 \[
 \xymatrix{C\ar[d]\ \ar@{^{(}->}[r]& {\C}\ar[d]_\alpha\ \ar@{^{(}->}[r]\ar[d] & {\Y}\ar[dl]^\beta \\
\mathrm{Spec}(\mathbb{C})\ar[r]^-{\{b\}} & B }
\]
 such that  the left square is   a smooth
 family of deformations of $C$ having surjective Kodaira-Spencer map at $p$; $\beta$ is a smooth family of
 projective surfaces and the upper right inclusion restricts over $b$ to the inclusion
 $C \subset Y$; moreover  $|\C(u)|$ is a positive-dimensional
 linear system on the surface $\Y(u)$ for all closed points $u \in B$.

 Let  $\L = \O_\Y(\C)$.  After possibly performing an etale base change, we
      can find a trivial free subsheaf of rank two
 $\O_B^{\oplus 2} \subset \beta_*\L$ which defines a rational $B$-map
 \begin{equation}\label{E:free4}
 \xymatrix{
{\Y}\ar@{-->}[r]^-\chi \ar[d]^\beta& {\P^1\times B}\ar[dl]\\
B}
 \end{equation}
  whose restriction over   $b$ is the rational map defined by the pencil  $\Lambda$.
  Let  $Z\subset \Y$ be the scheme of indeterminacy of $\chi$
  and let $\theta:\X \to \Y$ be the blow-up with center $Z$.
  Composing    with $\theta$ we obtain from  (\ref{E:free4}) a   family of deformations of $f$:
\[
\xymatrix{
X\ar[dd]^f \ \ar@{^{(}->}[r] & {\X}\ar[dd]^{\chi\cdot\theta} \ar[dr]^{\beta\cdot\theta} \\
&&B \\
{\P^1\times \{b\}} \ar@{^{(}->}[r] & {\P^1\times B}\ar[ur]}
\]
  After restricting over $(p,b)$ we obtain an induced family of deformations of $C$:
	\begin{equation}\label{E:free5}
	\xymatrix{
	C \ar[d]\ar[r] &\X \ar[d] \\
	\mathrm{Spec}(\mathbb{C}) \ar[r]^-{\{(p,b)\}}& \P^1\times B  }
	\end{equation}
	 Its Kodaira-Spencer map:
	\[
	\kappa(p,b): T_p\P^1\times T_bB \longrightarrow \mathrm{Ext}^1_C(\Omega^1_C,\O_C)
	\]
	is surjective because, by construction, it contains the image of the Kodaira-Spencer map of 
	$\alpha$. Now, recalling Lemma \ref{L:basech}, we deduce that Im$(\kappa(p,b))$ is contained in the image of the restriction 
	map
	\[
\xymatrix{
 H^0(\P^1,Ext^1_f(\Omega^1_{X/\P^1},\O_X))\ar[r]&
\mathrm{Ext}^1_C(\Omega^1_C,\O_C)}
\]
and therefore this map is surjective. Therefore
$Ext^1_f(\Omega^1_{X/\P^1},\O_X))$ is generated by its global
sections at $p$, thus $f$ is free.
  \qed

\begin{theorem}\label{T:unirfree}
  The following conditions are equivalent for an integer $g \ge 3$:
\begin{description}
\item[(i)] $\mgbar$ is   uniruled.

\item [(ii)] There exists a free fibration $f:X \longrightarrow
\P^1$
        with fibres of genus $g$.
				
		\item[(iii)] A general curve of genus $g$ moves in a
positive-dimensional linear system on some nonsingular projective surface which is not irrational ruled.
\end{description}
\end{theorem}

\proof  
$(i) \Longrightarrow (ii)$.  By assumption there is an algebraic
integral scheme $M$ of dimension $3g-4$ and a dominant rational
map
\[
\Psi:\xymatrix{ \P^1\times M \ar@{-->}[r] & \mgbar}
\]
 Since $\mgbar$ is projective we can take $M$   projective and
normalize it. Therefore $\Psi$ is defined on the complement of a
codimension two closed subset $Z$. Then  the image $\pi(Z)\subset
M$ under the projection $\pi: \P^1\times M \longrightarrow M$ is a
proper closed subset.  This means that, modulo replacing $M$ by an open subset, we
may assume that $\Psi$ is a morphism and that $M$ is nonsingular
of dimension $3g-4$. 
Let $(p,m)\in \P^1\times M$ be a general point. Then $\Psi(p,m) \in \overline{M}^\circ_g\subset \mgbar$,  
the open set of curves without automorphisms. the universal family over $\overline{M}^\circ_g$ pulls back to a family over an
open subset of $\P^1 = \P^1\times \{m\}$ containing $p$. After
embedding the total space into a projective surface and
desingularizing it we obtain a relatively minimal fibration $f:X\longrightarrow\P^1$
containing $C$ among its fibres.  Consider a semiuniversal
deformation of $f$:
\[
\xymatrix{
    X \ar@{^(->}[r]\ar[d]^f &\X \ar[d]^F \\
    \P^1 \ar@{^(->}[r]_-{\mathrm{id}\times \{v\}} & \P^1\times V }
    \]
    and the induced rational map
    \[
    \xymatrix{
    \psi_F: \P^1\times V \ar@{-->}[r] & \mgbar}
    \]
Since $\psi_F$ is well defined at the point $(p,v)$ and
$\psi_F(p,v)= \Psi(p,m)$ we see that $\psi_F$ is dominant because
$\Psi(p,m)$ is a general point of $\mgbar$.  It follows that the
differential of $\psi_F$ at $(p,v)$ is surjective. Recalling Lemma
\ref{L:basech}(iv) we deduce that the map:
\[
\xymatrix{
 r_p:H^0(\P^1,Ext^1_f(\Omega^1_{X/\P^1},\O_X))\ar[r]&
\mathrm{Ext}^1_{X(p)}(\Omega^1_{X(p)},\O_{X(p)})}
\]
is surjective.  But this means that
$Ext^1_f(\Omega^1_{X/\P^1},\O_X)$ is generated at $p$, and
therefore it is globally generated, i.e. $f$ is free.

\smallskip

$(ii) \Longrightarrow (iii)$.  By Lemma \ref{L:blowup} we may assume that $f$ is relatively minimal. By Theorem
\ref{T:semiuniv} it has a semiuniversal algebraic deformation
\[
    \xymatrix{
    X \ar@{^(->}[r]\ar[d]^f &\X \ar[d]^F \\
    \P^1 \ar@{^(->}[r]_-{\mathrm{id}\times \{v\}} & \P^1\times V }
\]
Since $f$ is free, $V$ is nonsingular at $v$. Let $p \in \P^1$ be such that the fibre $C=X(p)$ is nonsingular. Then we obtain  a family of deformations of $C$:
\[
\xymatrix{
C \ar@{^(->}[r]\ar[d] & \X \ar[d]^F \\
\mathrm{Spec}(\mathbb{C}) \ar[r]^-{\{p,v\}} & \P^1\times V }
\]
 and from Lemma \ref{L:basech} it follows that this family has surjective Kodaira-Spencer map at $v$.
Consider the composition 
\[
\xymatrix{
\X \ar[r]^-F &   \P^1\times V \ar[r] & V }
\]
and let $\Y := (\P^1\times V ) \times _V \X$. Then we have a commutative diagram:
\[                                                                  
\xymatrix{
\X \ar@{^(->}[r] \ar[dr]_-F& \Y\ar[d] \\
&\P^1\times V  }
\]
showing that $F$ defines   a family of general curves over $\P^1\times V$   moving in a linear pencil on a surface. 

\smallskip

$(iii) \Longrightarrow (i)$.  With the same notations as in the proof of Proposition \ref{C:free1}, the hypothesis 
implies that, for a general curve $C$ of genus $g$, there  a family of deformations (\ref{E:free5}). 
The induced functorial rational map
$\xymatrix{\P^1\times B \ar@{-->}[r]& \mgbar}$ is not constant along $\P^1\times \{b\}$ because 
the fibration $f$ is not isotrivial. Therefore $\mgbar$ is uniruled.
 \qed

   The
equivalence of conditions (i) and (iii) of the Theorem is  well-known (see the
Proposition on p. 25 of \cite{HM82}).

\section{The main estimate} \label{S:estim}

We will need the following:

\begin{lemma}\label{L:def1}
Let $C$ be a general curve of genus $g \ge 3$ contained in a
projective nonsingular  surface $Y$ which is not irrational ruled, and such that
$$
r:=\dim(|C|) \ge 1
$$
Let $f: X \to \P^1$ be the (free) fibration defined by a general pencil
contained in $|C|$ and containing $C$ as a fibre. Then:
\[
 h^0(\P^1,Ext^1_f(\Omega^1_{X/\P^1},\O_X)) \ge 3(g-1)+r+1
\]
\end{lemma}

\proof   There is an $(r-1)$-dimensional family of pencils
contained in $|C|$ and containing the curve $C$ as a fibre. This
family of pencils defines a family of free fibrations which has
dimension $r-1+2=r+1$ and is effectively parametrized, i.e. it has
injective Kodaira-Spencer map,  at $f$. In fact the pencils
containing $C$ as a member are parametrized by  the
$(r-1)$-dimensional characteristic linear series, defined by the
image of the restriction:
\[
\xymatrix{ H^0(\O_Y(C)) \ar[r] & H^0(\O_C(C))}
\]
At most finitely many fibrations coming from different pencils can
coincide, because the automorphism groups of the fibres are
finite.  Moreover
  the image of
the Kodaira-Spencer map is contained in
\[H^0(\P^1,Ext^1_f(\Omega^1_{X/\P^1},\O_X)(-1))
\]
 by the exact sequence (\ref{E:basech}) and therefore each
 pencil gives rise to a two-dimensional family of
fibrations parametrized by   the projectivities of $\P^1$ fixing
$p$.
 Therefore from the  freeness of $f$ it follows that:
\[ h^0(\P^1,Ext^1_f(\Omega^1_{X/\P^1},\O_X)) \ge r+1+ 3(g-1)\]
 \qed

We can now prove the following  refinement of the estimate
 (\ref{E:free1}).

\begin{theorem}\label{T:def1}
Let $C$ be a general nonsingular curve of genus $g \ge 3$ moving
in a positive dimensional linear system in a projective
nonsingular surface $Y$ which is not irrational ruled. Then:
\begin{equation}\label{E:surf1}
      10\chi(\O_Y) - 2 K_Y^2 \ge 4(g-1) -C^2 - h^0(K_Y-C)
\end{equation}
If moreover $\mathrm{dim}(|C|) \ge 2$ or $h^1(\O_C(2C))=0$ then
$H^0(Y,K_Y-C)=0$.
\end{theorem}

\proof  Let $\Lambda \subset |C|$ be a general pencil and let
$f:X\to\P^1$ be the fibration
 defined by   $\Lambda$.  Then   we have:
\[
K_X^2 = K_Y^2 - C^2, \qquad \chi(\O_X) = \chi(O_Y)
\]
Proposition \ref{P:def1} gives:
\begin{equation}\label{E:extf1}
h^0(\P^1, Ext^1_f(\Omega^1_{X/\P^1},\O_X)) = 11\chi(\O_Y) - 2
K_Y^2 +2C^2- 2(g-1)
\end{equation}
because $h^1(\P^1, Ext^1_f(\Omega^1_{X/\P^1},\O_X))=0$, since $f$
is free. From Lemma (\ref{L:def1})   we obtain:
\begin{equation}\label{E:extf2}
3(g-1)+h^0(\O_Y(C)) \le h^0(\P^1, Ext^1_f(\Omega^1_{X/\P^1},\O_X))
\end{equation}
Putting (\ref{E:extf1}) and (\ref{E:extf2}) together we get:
\begin{equation}
\label{E:surf0} 11\chi(\O_Y) - 2 K_Y^2 +2C^2 \ge 5(g-1)+
h^0(Y,\O_Y(C))
\end{equation}
By applying the Riemann--Roch   theorem we have:
$$
\chi(\O_Y(C)) = \chi(\O_Y) - (g-1-C^2)
$$
and therefore:
\[
 h^0(\O_Y(C)) \ge \chi(\O_Y) - (g-1-C^2)- h^0(Y,K_Y-C)
\]
 Substituting in (\ref{E:surf0}) we get
(\ref{E:surf1}).

Assume that $\mathrm{dim}(|C|) \ge 2$. Then $h^0(C,\O_C(C)) \ge 2$
and therefore $H^1(\O_C(2C))= 0$, by Proposition
\ref{P:gencur}(iii).    If $H^0(Y,K_Y-C)\ne 0$ then
$\O_Y(2C)\subset \O_Y(K_Y+C)$, and therefore $\O_C(2C) \subset
\omega_C$, i.e. $H^1(\O_C(2C))\ne 0$. This is a contradiction.
\qed

Theorem (\ref{T:def1}) will be applied to show that, for certain
   algebraic surfaces $Y$, there is an upper
   bound on the genus $g$ of
   a general curve which moves in a positive dimensional
   linear system on $Y$.

\begin{remark}\rm\label{R:def3}
    Note  that if $h^0(T_Y)=0$ then  the left hand side of  (\ref{E:surf1}) is
 \[
 10 \chi(\O_Y) - 2 K_Y^2 =-\chi(T_Y)=h^1(T_Y)-h^2(T_Y)=: \mu(Y)
 \]
 the expected number of moduli of the surface $Y$.
\end{remark}

\begin{example}\label{Ex:BV}\rm
In \cite{BV05}, \S 3,  it is shown that a general curve $C$ of
genus 15 can be embedded as a non-degenerate nonsingular curve of
degree 19 in $\P^6$, lying on a
  nonsingular canonical surface
    $Y \subset \P^6$ which is a complete intersection of 4 quadrics, and that dim$(|C|)=2$ on $Y$.
      Therefore, by
Theorem   \ref{T:def1},
     inequality (\ref{E:surf1}) holds and $H^0(K_Y-C)=0$. Let's check.
    The relevant numbers are in this case:
    \[
    C^2 =9, \qquad K_Y^2 = 16, \qquad \chi(\O_Y) = 8
    \]
  We find:
  \[
  10\chi(\O_Y) - 2K_Y^2 = 48 > 47 = 4(g-1)-C^2
  \]
 and this shows that  (\ref{E:surf1}) is sharp.
   \end{example}

\section{General curves on surfaces of non-negative  Kodaira dimension} \label{S:surf1}

In this section we start   analyzing the case of general curves on
non-rational surfaces.
 The following elementary
lemma   will be useful.

\begin{lemma}\label{L:KC}
 Let $\sigma: Y  \to Z$  be  a birational morphism of projective
 nonsingular surfaces and let  $D_0 \subset Z$ be  an  irreducible and reduced curve.
 Factor $\sigma$   as a sequence of blow-ups:
 \begin{equation}\label{E:factor1}
 \xymatrix{
 Y=Z_n \ar[r]^-{\sigma_n}&Z_{n-1}\ar[r]^-{\sigma_{n-1}}&\cdots\ar[r]^-{\sigma_2}&Z_1\ar[r]^-{\sigma_1}&Z}
 \end{equation}
 and   denote by $D_i \subset Z_i$ the proper transform of $D_0$ under
 \[
 \pi_i:=\sigma_1\cdots\sigma_i: \xymatrix{Z_i \ar[r]& Z,} \qquad i=1,\dots,n
 \]
   Assume that the center of $\sigma_i$ is a   point of $D_{i-1}$ of multiplicity $\nu_{i-1}\ge 1$, for each
 $i=1,\dots, n$. Then
\[
K_YD_n = K_ZD_0+\sum_{i=1}^n\nu_{i-1}
\]
In particular, if the center of $\sigma_i$ is a singular  point of
$D_{i-1}$ for all $i$, then
 \[
 K_YD_n \ge  K_ZD_0+2n
 \]
 \end{lemma}

\proof  The last assertion is an obvious consequence of the first.
It suffices to prove the first assertion in the case $n=1$.   We
have $K_Y = \sigma^*K_Z+E$ where $E$ is the   exceptional curve.
We have:
 \[
K_YD_1 = (\sigma^*K_Z+E)(\sigma^*D_0-\nu_0E)= K_ZD_0+\nu_0
\]
\qed

\begin{theorem}\label{C:def2}
Let $Y$ be a  projective nonsingular  surface   with
$\kappa$-dim$(Y) \ge 0$, and
 let $C\subset Y$ be a general   nonsingular curve of   genus $g\ge 3$
 moving in a positive dimensional linear system.
 Then
\[g \le 5p_g(Y)+6+ {1\over 2}h^0(K_Y-C)\]
In particular:
\[
g \le \begin{cases}6,&\text{if  $p_g=0$} \\ 11,&\text{if $p_g=1$}
\end{cases}
\]
 \end{theorem}

 \proof   Let $\sigma: \xymatrix{Y\ar[r]&Z}$ be
 the birational morphism onto the minimal model of $Y$.
 After possibly contracting finitely many $(-1)$-curves
 on $Y$ we may assume that $\sigma$ can be factored as a sequence of $\delta \ge 0$ blow-ups
\[
 \xymatrix{
 Y=Z_\delta \ar[r]^-{\sigma_\delta}&Z_{\delta-1}\ar[r]^-{\sigma_{\delta-1}}&
 \cdots\ar[r]^-{\sigma_2}&Z_1\ar[r]^-{\sigma_1}&Z}
 \]
 so that, for each $i=1,\dots, \delta$, letting $D_i = (\sigma_{i+1}\cdots\sigma_\delta)(C) \subset Z_i$,
  the center of $\sigma_i:Z_i \to Z_{i-1}$ is a singular point of $D_{i-1}$.

 In all three cases inequality (\ref{E:surf1}) holds. From Lemma \ref{L:KC} it follows that
 \[
 C^2 = 2(g-1)-K_YC \le 2(g-1) - [\sigma(C)K_Z+2\delta] \le 2(g-1) - 2\delta
 \]
 because $\sigma(C)K_Z \ge 0$ by the assumption on the Kodaira dimension of $Z$.
 Therefore, since $K_Y^2= K_Z^2-\delta$, using (\ref{E:surf1}) we obtain:
 \[
 \begin{array}{lll}
 10\chi(\O_Z)-2K_Z^2+2\delta &= 10\chi(\O_Y) - 2K_Y^2 \\
 & \ge 4(g-1)-C^2 -h^0(K_Y-C)  \\
 &\ge 4(g-1) -[2(g-1)-2\delta]-h^0(K_Y-C)
  \end{array}
 \]
 This  gives
 \[
 10\chi(\O_Z)-2K_Z^2 \ge 2(g-1)-h^0(K_Y-C)
 \]
 The conclusion now follows  because $\chi(\O_Z)\le 1+p_g$ and $K_Z^2 \ge 0$.  \qed

  \begin{theorem}\label{T:ell1}
Let $Y$ be a non-ruled and non-rational   elliptic surface, $\pi:Y
\to B$ the elliptic fibration onto a nonsingular connected curve
$B$.  Assume that there exists    a general nonsingular connected
curve  $C \subset Y$  of  genus $g\ge 3$ such that dim$(|C|) \ge
1$.  Then $B=\P^1$ and  $g \le 16$.
\end{theorem}

\proof
 Let $|K| = |M|+D$, where $|M|$ is the mobile part and $D$ is the fixed divisor.
Then $M= \sum_{i=1}^m F_{b_i}$,  where $b_1,\dots,b_m\in B$ and
$F_{b_i}=\pi^{-1}(b_i)$. It follows that $h^0(K_Y-C)=0$ because
all fibres of $\pi$ have arithmetic genus $1$. 

Let $F$ be a general fibre of $\pi$. We have $C\cdot F =: k \ge
1$.  If $k=1$ then $\pi_{|C}: C \to B$ is an isomorphism, and it
follows that every linear pencil $\Lambda$ containing $C$ is
isotrivial. Since $C$ is general, it follows that $Y$ is ruled
(Proposition   \ref{P:free1}), which is a contradiction.

Therefore $k \ge 2$.   If $B$ has genus $g(B) \ge 1$ then
 \[
 \pi_{|C}: \xymatrix{ C \ar[r]& B}
 \]
  is an irrational involution on
 $C$, and this contradicts the fact that $C$ has general moduli (Proposition \ref{P:gencur}(iv)).
Therefore $B \cong \P^1$, and $ \pi_{|C}$ defines a $g^1_k$ on
$C$;
 in other words  $h^0(C,\O_C(F)) \ge 2$. 
From Theorem
\ref{C:def2}  we deduce that
  $g  \le 16$ if $p_g \le 2$.
Therefore we may assume that $p_g \ge 3$.
Since $p_g \ge 3$,
  we have $m \ge 2$, and it follows that
  \[
  \O_C(mF) \subset \O_C(K_Y) \subset
   \O_C(K_Y+C) = \omega_C
  \]
   This means that $h^1(C,\O_C(mF))\ne 0$, while
   $h^0(C,\O_C(F)) \ge 2$,     contradicting the generality of
  $C$ (Proposition \ref{P:gencur}(iii)).    \qed

\begin{remarks}\rm\label{R:bound}
  The  well-known bound $g \le 11$ for  the genus of a \emph{nonsingular}
  curve  with general moduli on a K3 surface
  can be  deduced from a dimension  count,  as shown  in \cite{aM72}.
   In \cite{MM83}  and
  \cite{sM88} it is proved that a general nonsingular curve of genus $g$
  can be embedded in a K3 surface if and only if  $g \le 11$  and $g\ne 10$.
   In particular the bound
 of Theorem \ref{C:def2} is   sharp for the case $p_g=1$.
 In \cite{FKPS07} it is shown that even a general curve of genus 10 can
 be    birationally embedded in a K3 surface
 provided one allows it to have one node (hence arithmetic genus 11).
 I don't know if the bound $g \le 16$ in Theorems \ref{C:def2} and
\ref{T:ell1}  is sharp.
  \end{remarks}

  The following result takes care of  irregular surfaces of positive geometric
  genus, showing the impossibility for a general curve of genus $g \ge 3$ to move
  in a positive dimensional linear system on such a surface.
  The case $p_g=0$ is already included in previous results
  of this section.

  \begin{theorem}\label{T:irregular}
Let $Y$ be a projective nonsingular surface with  $p_g>0$ and
$q>0$, and  let $C\subset Y$ be a nonsingular curve of genus $g
\ge 3$ such that  $\dim(|C|) \ge 1$.  Then no pencil $\Lambda
\subset |C|$ containing $C$ as a member defines a free fibration.
\end{theorem}

\proof Let  $f:X \longrightarrow \P^1$ be the fibration associated
to a pencil $\Lambda$ as in the statement. We have:
\[
h^2(T_X) = h^2(T_Y) \ge p_g+q-1 \ge p_g
\]
The equality is a standard computation (see e.g. \cite{eS06},
Example 3.4.13(iv)), and the first inequality is a classical
result of B. Segre \cite{bS34} (see \cite{dM71}, p. 127 for a
modern proof). Now consider the exact sequence (\ref{E:reldiff3})
associated to $f$. If  $f$ is free then
$H^1(\P^1,Ext^1_f(\Omega^1_{X/\P^1},\O_X))= 0$, and therefore
$H^2(T_X) \cong H^1(T_{\P^1}\otimes R^1f_*\O_X)$.  But
\[
p_g = h^1(R^1f_*\O_X) > h^1(T_{\P^1}\otimes
R^1f_*\O_X) = h^2(T_X)
\]
Here the inequality is a consequence of the assumption $p_g>0$ and the last equality follows from the exact sequence 
(\ref{E:reldiff3}). Therefore  we get a contradiction. \qed

  \section{General curves on surfaces of general type} \label{S:surfgt}

  The following is our main result on general curves on surfaces
  of general type.

   \begin{theorem}\label{T:gt1}
Let  $Y$ be a   surface of general type and let  $Z$ be the
minimal model of $Y$. Assume that $K_Z^2 \ge 3\chi(\O_Z) - 10$ and
that $C \subset Y$  is a general nonsingular  connected curve of
genus $g \ge 3$.   If     one of the following holds:
\begin{itemize}
\item[(a)] $\dim(|C|) \ge 2$,

\item[(b)]  $\dim(|C|)=1$, $h^0(K_Y-C)=0$ and  $C^2\ge {(g-1)\over
2}$,

\item[(c)] $\dim(|C|)=1$, and $h^1(C,\O_C(2C))=0$.
\end{itemize}
Then
 $g \le 19$.
 \end{theorem}

\proof   By Theorem \ref{C:def2} we may assume that $Y$ has
geometric genus $p_g \ge 2$. As in the proof of Theorem
\ref{C:def2} we may assume that the birational morphism
\[
\sigma: \xymatrix{Y \ar[r] & Z}
\]
can be factored as a sequence of $\delta \ge 0$ blow-ups
\[
 \xymatrix{
 Y=Z_\delta \ar[r]^-{\sigma_\delta}&Z_{\delta-1}\ar[r]^-{\sigma_{\delta-1}}
 &\cdots\ar[r]^-{\sigma_2}&Z_1\ar[r]^-{\sigma_1}&Z}
 \]
 so that, for each $i=1,\dots, \delta$, letting $D_i = (\sigma_{i+1}\cdots\sigma_\delta)(C) \subset Z_i$,
  the center of $\sigma_i:Z_i \to Z_{i-1}$ is a singular point of $D_{i-1}$.  We have
 $K_Y = \sigma^*K_Z(\sum r_iE_i)$, where the $E_i$'s are the
  irreducible components of the exceptional locus of $\sigma$ and $r_i > 0$ for all $i$.

Let $L =\O_C(\sigma^*K_Z)$.  Then, since   $H^0(K_Y-C)=0$ (by
Theorem \ref{T:def1})  we have $H^0(\sigma^*K_Z-C) =0$ and
therefore
 \begin{equation}\label{E:gt2}
 h^0(L) \ge h^0(\sigma^*K_Z) = p_g \ge 2
 \end{equation}
Inequality (\ref{E:surf1}) can be written under the form:
\[
10\chi(\O_Z)-2K_Z^2+2\delta \ge 4(g-1)- C^2
\]
After substitution of the inequality  $K_Z^2 \ge 3\chi(\O_Z) - 10$
we obtain:
\begin{equation}\label{E:gt4}
4\chi(\O_Z)+20+2\delta \ge 4(g-1)-C^2
\end{equation}
which implies:
\[
28+2\delta \ge 4(g-1) - C^2 - 4p_g +4
\]
From Lemma \ref{L:KC} we deduce
\[
2(g-1)-C^2 = \deg(\O_C(K_Y)) \ge \deg(L)+2\delta
\]
and therefore we obtain:
\begin{equation}\label{E:gt0}
    28+2\delta  \ge 2[\deg(L)-2p_g+2+2\delta]  + C^2
\end{equation}
Assume that we are in case a). Then
\begin{equation}\label{E:gt1}
h^0(\omega_CL^{-1}) = h^0(\O_C(C+\sum r_iE_i)) \ge h^0(\O_C(C))
\ge \dim(|C|) \ge 2
\end{equation}
Inequalities (\ref{E:gt1}) and (\ref{E:gt2}) imply that $L$
contributes to the Clifford index of $C$.
 Since $C$ is general,  recalling that
\[
{\rm Cliff}(L) = \deg(L) - 2h^0(L)+2
\]
we have the following relations satisfied by the  Clifford
indices:
 \begin{equation}\label{E:gt3}
\bigg[{g-1 \over 2}\bigg] ={\rm Cliff}(C)  \le {\rm Cliff}(L) \le
\deg(L) -2p_g+2
\end{equation}
For the   equality see Proposition \ref{P:gencur}(i).   The first
inequality follows from the very definition of Cliff$(C)$, and the
last   is by (\ref{E:gt2}).
 Substituting in (\ref{E:gt0}) we obtain:
\[
\begin{array}{llll}
28+2\delta  &\ge 2[\deg(L)-2p_g+2+2\delta]  + C^2 \\ \\
 & \ge 2[{\rm Cliff}(C) +2\delta]+C^2 \\ \\
& \ge g-2 + 4\delta+C^2
\end{array}
\]
Since $h^0(\O_C(C)) \ge 2$, from Proposition \ref{P:gencur}(ii) we
have $C^2 \ge {1\over 2}g + 1$.  Therefore:
\[
28 -2\delta  \ge {3\over 2} g -1
\]
which implies $g \le 19$.

Assume now that we are in case b).  We may assume that
$h^0(\omega_CL^{-1})=1$, because otherwise $L$ contributes to
$\mathrm{Cliff}(C)$ and we conclude as before.

Then:
\[
\deg(L)-2p_g+2 \ge \deg(L)-2h^0(L)+2 = \deg(\omega_CL^{-1})\ge C^2
\]
Substituting in (\ref{E:gt0}) we obtain:
\[
\begin{array}{llll}
28+2\delta  &\ge 2[\deg(L)-2p_g+2+2\delta]  + C^2 \\ \\
 & \ge 2(C^2+2\delta)+C^2 = 3C^2+4\delta \ge {3\over 2}(g-1)+4\delta
\end{array}
\]
because of the assumptions.  Therefore we obtain  $g \le 19$ again.

If we are in case c) then $C^2 \ge g/2$, since $\O_C(2C)$ is
non-special, and $H^0(K_Y-C)=0$ by Theorem \ref{T:def1}. Therefore
we reduce to case b).   \qed


\noindent\textsc{address of the author:}

\noindent Dipartimento di Matematica e Fisica,
  Universit\`a Roma Tre \newline Largo S. L. Murialdo 1,
  00146 Roma, Italy.
  \smallskip
\newline  \texttt{sernesi@mat.uniroma3.it}

\end{document}

%% file: amsmacros.tex
\def\P{\mathbb{P}}

\def\O{\mathcal{O}}
\def\k{{\bf k}}

\def\X{\mathcal{X}}

\def\Y{\mathcal{Y}}

\def\L{\mathcal{L}}

\def\C{\mathcal{C}}
\def\M{\mathcal{M}}

\def\dim{\mathrm{dim}}

\def\separation{\medskip}
  \newtheorem{theorem}{Theorem}[section]
\newtheorem{lemma}[theorem]{Lemma}
\newtheorem{prop}[theorem]{Proposition}
\newtheorem{definition}[theorem]{Definition}

\newtheorem{remarks}[theorem]{Remarks}
\newtheorem{remark}[theorem]{Remark}

\newtheorem{example}[theorem]{Example}

\newcommand{\proof}{{\it Proof.}\ }
\newcommand{\qed}{\hfill  $\Box$\separation}